\documentclass[a4paper,12pt]{article}
\usepackage{fullpage}
\usepackage[english,francais]{babel}
\usepackage{epsf}
\usepackage[dvips]{epsfig}
\usepackage{amssymb}
\usepackage{amsmath}
\usepackage[all]{xy}

\usepackage{ucs}
\usepackage[utf8x]{inputenc}\usepackage{wasysym}
\usepackage{aeguill}

\newtheorem{theorem}{Theorem}
\newtheorem{prop-f}[theoreme]{Proposition}
\newtheorem{prop}[theorem]{Proposition}

\newtheorem{lemma}[theorem]{Lemma}

\newcommand{\finpreuve}{\hspace{\stretch{1}}{$\square$}}

\newcommand{\R}{\mathbb{R}}

\def\esp{\medskip\noindent}
\def\prg#1{\esp{\bf #1. }}

\def\proof{\prg{Proof}}

\def\proofof#1{\prg{Proof of #1}}
\renewcommand{\epsilon}{\varepsilon}
\def\btab{\begin{eqnarray*}}
\def\etab{\end{eqnarray*}}
\def\beq{\begin{equation}}
\def\eeq{\end{equation}}

\newcounter{numeroexo}

\def\perco#1#2#3{#1 \; \leftrightarrow_{#3} \; #2}
\def\percorayon#1#2#3#4{#1 \;\leftrightarrow_{#3}^{\le #4} \; #2}

\begin{document}

\selectlanguage{english}
\title{Percolation in a multiscale Boolean model}

\author{
Jean-Baptiste Gou\'er\'e
\footnote{
\noindent\textit{Postal address}:
Universit\'e d'Orl\'eans
MAPMO
B.P. 6759
45067 Orl\'eans Cedex 2
France
\textit{E-mail}: jbgouere@univ-orleans.fr
}
}

\date{}

\maketitle

\begin{abstract}
We consider percolation in a multiscale Boolean model.
This model is defined as the union of scaled independent copies of a given Boolean model.
The scale factor of the $n^{\textrm{th}}$ copy is $\rho^{-n}$.
We prove, under optimal integrability assumptions, that no percolation occurs in the multiscale Boolean model for large enough $\rho$ 
if the rate of the Boolean model is below some critical value.
\end{abstract}

\section{Introduction and statement of the main result}

\subsection{The Boolean model}

Let $d\ge 2$.
Let $\mu$ be a finite measure on $]0,+\infty[$.
We assume that the mass of $\mu$ is positive. 
Let $\xi$ be a Poisson point process on $\R^d \times ]0,+\infty[$ whose intensity is the product of the Lebesgue measure on $\R^d$ by $\mu$.
With $\xi$ we associate a random set $\Sigma(\mu)$ defined as follows:
$$
\Sigma(\mu)=\bigcup_{(c,r) \in \xi} B(c,r)
$$
where $B(c,r)$ is the open Euclidean ball of radius $r$ centered at $c$.
The random set $\Sigma(\mu)$ is the Boolean model with parameter $\mu$.
When shall sometimes write $\Sigma$ to simplify the notations.

The following description may be more intuitive.
Let $\chi$ denote the projection of $\xi$ on $\R^d$.
With probability one this projection is one-to-one.
We can therefore write:
$$
\xi=\{(c,r(c)), c \in \chi\}.
$$
Write $\mu=m\nu$ where $\nu$ is a probability measure.
Then, $\chi$ is a Poisson point process on $\R^d$ with density $m$.
Moreover, given $\chi$, the sequence $(r(c))_{c\in\chi}$ is a sequence of independent random variable with common distribution $\nu$.
We shall not use this point of view.

\subsection{Percolation in the Boolean model}
\label{s:boolean}

Let $C$ denote the connected component of $\Sigma$ that contains the origin.
We say that $\Sigma$ percolates if $C$ is unbounded with positive probability.
We refer to the book by Meester and Roy  \cite{Meester-Roy-livre} for background on continuum percolation.
Set:
$$
\lambda_c(\mu) = \inf \{\lambda>0 : \Sigma(\lambda\mu) \hbox{ percolates}\}.
$$
One easily check that $\lambda_c(\mu)$ is finite as soon as $\mu$ has a positive mass.
In \cite{G-perco-boolean-model} we proved that $\lambda_c(\mu)$ is positive if and only if: 
$$
\int r^d\mu(dr) < \infty.
$$
The only if part had been proved earlier by Hall \cite{Hall-continuum-percolation}. 
For all $A,B \subset \R^d$, we write $\perco{A}{B}{\Sigma}$  if there exists a path in $\Sigma$ from $A$ to $B$.
We denote by $S(c,r)$ the Euclidean sphere or radius $r$ centered at $c$ :
$$
S(c,r)=\{x \in \R^d : \|x-c\|_2=r\}.
$$
We write $S(r)$ when $c=0$.

The critical parameter $\lambda_c(\mu)$ can also be defined as follows:
$$
\lambda_c(\mu) 
   =  \sup \left\{ \lambda>0 : P\left(\perco{\{0\}}{S(r)}{\Sigma(\lambda\mu)}\right) \to 0 \hbox{ as } r \to \infty\right\}, \\
$$
We shall need two other critical parameters:
\begin{eqnarray*}
\widehat{\lambda_c}(\mu) 
  & = & \sup \left\{ \lambda>0 : P\left(\perco{S(r/2)}{S(r)}{\Sigma(\lambda\mu)}\right) \to 0 \hbox{ as } r \to \infty\right\}, \\
\widetilde{\lambda_c}(\mu) 
  & = & \sup \left\{ \lambda>0 : r^dP\left(\perco{\{0\}}{S(r)}{\Sigma(\lambda\mu)}\right) \to 0 \hbox{ as } r \to \infty\right\}.
\end{eqnarray*}
We have (see Lemma \ref{l:lambdacinegalite}) :
\begin{equation}\label{e:lambdacinegalite}
\widetilde{\lambda_c}(\mu) \le \widehat{\lambda_c}(\mu) \le \lambda_c(\mu). 
\end{equation}

When the support of $\mu$ is bounded, 
$$
P\left(\perco{\{0\}}{S(r)}{\Sigma(\lambda\mu)}\right)
$$
decays exponentially fast to $0$ as soon as $\lambda<\lambda_c(\mu)$
(see for example \cite{Meester-Roy-livre}, Section 12.10 in \cite{Grimmett-percolation} in the case of constant radii
or the papers \cite{Meester-Roy-Sarkar}, \cite{Menshikov-Sidorenko-coincidence}, \cite{Zuev-I} and \cite{Zuev-II}).
Therefore:
\begin{equation}\label{e:lambdacborne}
\widetilde{\lambda_c}(\mu)=\widehat{\lambda_c}(\mu)=\lambda_c(\mu) \hbox{ as soon as the support of }\mu\hbox{ is bounded}.
\end{equation}

\prg{Remarks} 
\begin{itemize}
\item
The treshold parameter $\widehat{\lambda_c}(\mu)$ is positive if and only if $\int r^d \mu(dr)$ is finite 
(i.e., if and only if $\lambda_c(\mu)$ is positive). See Lemma \ref{l:hat}.
\item Using ideas of \cite{G-perco-boolean-model},
we can check that $\widetilde{\lambda_c}(\mu)$ is positive if and only if 
$$
x^d \int_x^\infty r^d \mu(dr) \to 0 \hbox{ as } x\to\infty.
$$
If we only use results stated in \cite{G-perco-boolean-model}, we can easily get the following weaker statements.
Let $D(\lambda\mu)$ denote the Euclidean diameter of the connected component of $\Sigma(\lambda\mu)$ that contains the origin.
Note that $\widetilde{\lambda_c}(\mu)$ is positive if and only if there exists $\lambda$ such that:
\begin{equation}\label{e:decay}
r^dP(D(\lambda\mu) \ge r) \to 0, \hbox{ as } r\to\infty.
\end{equation}
If $E(D(\lambda\mu)^d)$ is finite then \eqref{e:decay} holds.
If \eqref{e:decay} holds then $E(D(\lambda\mu)^{d-\epsilon})$ is finite for any small enough $\epsilon>0$.
By Theorem 2.2 of \cite{G-perco-boolean-model} we thus get the following implications:
$$
\int_0^{+\infty} r^{2d} \mu(dr)<\infty \; \hbox{ implies } \; \widetilde{\lambda_c}(\mu)>0 
\; \hbox{ implies }\; \forall\epsilon>0: \int_0^{+\infty} r^{2d-\epsilon} \mu(dr)<\infty.
$$

\end{itemize}

\subsection{A multiscale Boolean model}

Let $\rho>1$ be a scale factor. 
Let ${(\Sigma_n)}_{n\ge 0}$ be a sequence of independent copies of $\Sigma(\mu)$.
In this paper, we are interested in percolation properties of the following multiscale Boolean model:
\begin{equation}\label{e:definition-multi}
\Sigma^{\rho}(\mu)=\bigcup_{n \ge 0} \rho^{-n} \Sigma_n.
\end{equation}
We shall sometimes write $\Sigma^{\rho}$ to simplify the notations.
As before,
we say that $\Sigma^{\rho}$ percolates if the connected component of $\Sigma^{\rho}$ that contains the origin is unbounded with positive probability.

This model seems to have been first introduced as a model of failure in geophysical medias in the $80'$.
We refer to the paper by Molchanov, Pisarenko and Reznikova \cite{Molchanov-al-failure} for an account of those studies.
For more recent results we refer to 
\cite{Broman-Camia-self-similar}, \cite{Meester-Roy-livre}, \cite{Meester-Roy-Sarkar}, \cite{Menshikov-al-multi}, 
\cite{Menshikov-al-multi-unbounded} and \cite{Popov-V-sticks}.

This model is related to a discrete model introduced by Mandelbrot \cite{Mandelbrot-fractale-perco}.
We refer to the survey by L. Chayes \cite{Chayes-aspects-fractale-perco} 
and, for more recent results, to \cite{Broman-Camia-mandelbrot}, \cite{O-fractale-perco} and \cite{White}.

In \cite{Menshikov-al-multi}, Menshikov, Popov and Vachkovskaia considered the case where the radii of the unscaled process $\Sigma_0$ equal $1$.
They proved the following result.
\begin{theorem}[\cite{Menshikov-al-multi}] \label{th:MPV1}
If $\lambda<\lambda_c(\delta_1)$ then, for all large enough $\rho$, $\Sigma^{\rho}(\lambda\mu)$ does not percolate.
\end{theorem}

In \cite{Menshikov-al-multi-unbounded} the same authors considered the case where the radii are random and can be unbounded.
They considered the following sub-autosimilarity assumption on the measure $\mu$:
\begin{equation}\label{e:MPV2}
\lim_{a\to\infty} \sup_{r\ge 1/2} \frac{a^d\mu([ar,+\infty[)}{\mu([r,+\infty[)}=0
\end{equation}
with the convention $0/0=0$.
They proved the following result.
\begin{theorem}[\cite{Menshikov-al-multi-unbounded}] \label{th:MPV2}
Assume that the measure $\mu$ satisfies \eqref{e:MPV2}.
Assume that $\widetilde{\lambda_c}(\mu)$ is positive. 
If $\lambda < \widetilde{\lambda_c}(\mu)$ then, for all large enough $\rho$, $\Sigma^{\rho}(\lambda\mu)$ does not percolate.
\end{theorem}

Note that \eqref{e:MPV2} is fulfilled for any measure with bounded support.
Because of \eqref{e:lambdacborne}, Theorem \ref{th:MPV2} is then a generalization of Theorem \ref{th:MPV1}.

In \cite{G-perco-generale} we proved the following related result in which $\rho$ is fixed.

\begin{theorem}[\cite{G-perco-generale}] \label{th:Gmulti}
Let $\rho>1$.
There exists $\lambda>0$ such that $\Sigma^{\rho}(\lambda\mu)$ does not percolate if and only if:
\begin{equation}\label{e:th}
\int_{[1,+\infty[} \beta^d \ln(\beta) \mu(d\beta)<\infty.
\end{equation}
\end{theorem}

The main result of this paper is the first item of the following theorem.
The second item is easy and already contained in Theorem \ref{th:Gmulti}.
Recall that, by Lemma \ref{l:hat}, $\widehat{\lambda_c}(\mu)$ is positive as soon as $\int r^d\mu(dr)$ is finite and therefore as soon as 
\eqref{e:th} holds.
\begin{theorem} \label{th:0}
$ $
\begin{enumerate}
\item Assume \eqref{e:th}.
Then, for all $\lambda<\widehat{\lambda}_c(\mu)$, there exists $\rho(\lambda)>1$ such that, for all $\rho\ge\rho(\lambda)$:
\begin{equation}\label{e:pitilde0}
P\left(\perco{S(r/2)}{S(r)}{\Sigma^{\rho}(\lambda\mu)}\right) \to 0 \hbox{ as } r \to \infty
\end{equation}
and therefore $\Sigma^{\rho}(\lambda\mu)$ does not percolate.
\item 
Assume that \eqref{e:th} does not hold.
Then, for all $\lambda>0$ and for all $\rho>1$, $\Sigma^{\rho}(\lambda\mu)$ percolates.
\end{enumerate}
\end{theorem}

The proof is given in Section \ref{s:preuve-th:0}. 
The ideas of its proof and the ideas of the proofs of Theorems \ref{th:MPV1} and \ref{th:MPV2} are given in Subsection \ref{s:idees}.

The first item of Theorem \ref{th:0} is a generalization of Theorem \ref{th:MPV2} and thus of Theorem \ref{th:MPV1}.
Indeed, by \eqref{e:lambdacinegalite}, one has $\lambda<\widehat{\lambda_c}$ as soon as $\lambda<\widetilde{\lambda_c}$.
Moreover, by the second item of Theorem \ref{th:0}, \eqref{e:th} has to be a consequence of the assumptions of Theorem \ref{th:MPV2}.
For example, one can check that \eqref{e:th} is a consequence of \eqref{e:MPV2}
\footnote{From \eqref{e:MPV2} one gets the existence of $a>1$ such that, for all $r \ge a$, one has $\mu([r,+\infty[)\le 2^{-1}a^{-d}\mu([r/a,+\infty[)$.
By induction and standard computations this yields, for all $r \ge a$, $\mu([r,+\infty[) \le Ar^{-\ln(2)/\ln(a)-d}$.
Therefore, for a small enough $\eta>0$, one has $\int r^{d+\eta} \mu(dr)<\infty$.}.
Alternatively, one can check that \eqref{e:th} is a consequence of $\widetilde{\lambda_c}(\mu)>0$ (see the remarks at the end of Section \ref{s:boolean}).

Let us denote by $\lambda_c(m^{\rho}_{\infty})$ and $\widehat{\lambda}_c(m^{\rho}_{\infty})$ the $\lambda_c$ and $\widehat{\lambda}_c$
critical tresholds for the multiscale model with scale parameter $\rho$.
Theorems \ref{th:Gmulti} and \ref{th:0} yield the following result:
\begin{enumerate}
\item If \eqref{e:th} holds then $\lambda_c(m^{\rho}_{\infty})>0$ (and actually the proof of Theorem \ref{th:Gmulti} yields
$\widehat{\lambda}_c(m^\rho_{\infty})>0$) for all $\rho>1$ and
$\widehat{\lambda}_c(m^{\rho}_{\infty}) \to \widehat{\lambda}_c(\mu) > 0$ as $\rho \to \infty$.
\item Otherwise, $\widehat{\lambda}_c(m^{\rho}_{\infty})=\lambda_c(m^{\rho}_{\infty}) = 0$ for all $\rho>1$.
\end{enumerate}

Let us denote by $D^{\rho}(\lambda\mu)$ the diameter of the connected component of $\Sigma^{\rho}(\lambda\mu)$ that contains the origin.
The following result is an easy consequence of Theorem \ref{th:0} above and Theorems 2.9 and 1.2 in \cite{G-perco-generale}.

\begin{theorem} \label{th:1}
Let $s>0$, $\lambda>0$ and $\rho>1$.
\begin{enumerate}
\item If $\int_{[1,+\infty[} \beta^{d+s} \mu(d\beta)<\infty$ and \eqref{e:pitilde0} holds, then $E\big((D^{\rho}(\lambda\mu))^s\big)<\infty$.
\item If $\int_{[1,+\infty[} \beta^{d+s} \mu(d\beta)=\infty$ then $E\big((D^{\rho}(\lambda\mu))^s\big)=\infty$.
\end{enumerate}
\end{theorem}

The proof is given is Section \ref{s:preuve-th:1}.

\subsection{Superposition of Boolean models with different laws}

Using the same arguments as in the proof of Theorem \ref{th:0}, we could prove similar results for infinite superpositions 
$$
\bigcup_{n \ge 0} \rho^{-n} \Sigma_n
$$ 
where the Boolean models $\Sigma_n$ are independent but not identically distributed.
We will not give such a result here.
However, we wish to give a weaker result for the superposition of two independent Boolean models at different scales.
As we consider only two scales the proof is easier than the proof of Theorem \ref{th:0}.
The proof uses Lemmas \ref{l:eta}, \ref{l:dichotomie} and \ref{l:carre} and is given in Section \ref{s:preuve-p}.
The result gives some insight on the critical treshold in the case of balls of random radii.

This result, in the case where the supports of $\nu_1$ and $\nu_2$ are bounded, is already implicit in \cite{Meester-Roy-Sarkar}
in their proof of non universality of critical covered volume (see \eqref{e:phic} below).
See also \cite{Molchanov-al-failure}.

\begin{prop} \label{p}
Let $\nu_1$ and $\nu_2$ be two finite measures on $]0,+\infty[$.
We assume that the masses of $\nu_1$ and $\nu_2$ are positive.
Let $0<\alpha<1$.
Then, for all $\rho>1$,
$$
\widehat{\lambda}_c(\alpha \nu_1 + (1-\alpha) H_{\rho} \nu_2) 
 \le
\min\big( \widehat{\lambda}_c(\alpha \nu_1),\widehat{\lambda}_c((1-\alpha) H_{\rho} \nu_2)\big)
 =
\min\left(\frac{\widehat{\lambda}_c(\nu_1)}{\alpha}, \frac{\widehat{\lambda}_c(\nu_2)}{1-\alpha}\right).
$$
Moreover,
$$
\widehat{\lambda}_c(\alpha \nu_1 + (1-\alpha) H_{\rho} \nu_2) 
 \to 
\min\left(\frac{\widehat{\lambda}_c(\nu_1)}{\alpha}, \frac{\widehat{\lambda}_c(\nu_2)}{1-\alpha}\right)
\hbox{ as } \rho\to\infty.
$$
The above convergence is uniform in $\alpha$.
\end{prop}

\begin{figure}[h!] 
\centering
\includegraphics[width=0.95\textwidth]{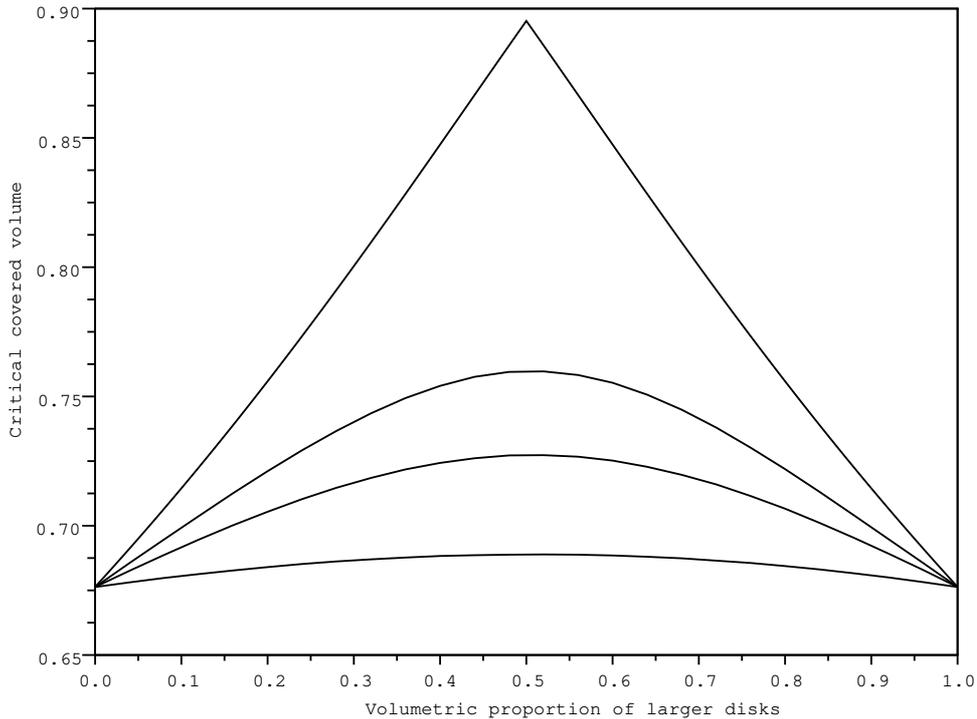}
\caption{Critical covered volume as a function of $\alpha$ for different values of $\rho$. 
From bottom to top: $\rho=2, \rho=5, \rho=10$ and the limit as $\rho \to \infty$.}
\label{f}
\end{figure}

We now make some remarks about this result and about some related numerical results.
For a finite measure $\mu$ on $]0,+\infty[$, we denote by $\phi_c(\mu)$ the critical covered volume:
\begin{equation}\label{e:phic}
\phi_c(\mu) = P\big( 0 \in \Sigma(\lambda_c(\mu)\mu)\big) = 1-\exp\left(-\lambda_c(\mu)\int v_dr^d\mu(dr)\right)
\end{equation}
where $v_d$ is the volume of the unit Euclidean ball in $\R^d$.
This is the mean volume occupied by the critical Boolean model and this is scale invariant.
Let us assume that $\nu_1=\nu_2=\delta_1$.
By \eqref{e:lambdacborne}, by Proposition \ref{p} and with the above notation we have:
\begin{equation}\label{e:fatigue}
\phi_c(\alpha \delta_1 + (1-\alpha) H_{\rho} \delta_1) 
 \to 
1-\exp\left(-v_d \lambda_c(\delta_1)\min\left(\frac{1}{\alpha}, \frac{1}{1-\alpha}\right)\right).
\end{equation}
There are several numerical studies of the above critical covered volume when $d=2$ and $d=3$.
To the best of our knowledge, the most acccurate values when $d=2$ are given in \cite{QZ-PRE-2007}.
Let us assume henceforth that $d=2$.
In \cite{QZ-PRE-2007}, the authors give: 
\begin{equation}\label{e:approx}
\phi_c(\delta_1)=1-\exp(-v_2 \lambda_c(\delta_1)) \approx 0.6763475(6).
\end{equation}
In Figure \ref{f} we reproduce the graph of critical covered volume $\phi(\alpha,\rho)$ as a function of $\alpha$ when $\rho=2$, $\rho=5$ and $\rho=10$
(see \cite{QZ-PRE-2007} for more results).
We also represent the graph of the right-hand side of \eqref{e:fatigue}, that we denote by $\phi(\alpha,\infty)$, as a function of $\alpha$.
We use \eqref{e:approx} to get an approximate value of $v_2\lambda_c(\delta_1)$.

\paragraph{Remarks} 
\begin{itemize}
\item 
When $\rho \to \infty$, 
the critical covered volume $\phi(\cdot,\rho)$ converges to $\phi(\cdot,\infty)$ which is symmetric: $\phi(\alpha,\infty)=\phi(1-\alpha,\infty)$.
When $\rho$ is finite, the critical covered volume may also look symmetric but Quintanilla and Ziff showed, 
based on their numerical simulations and statistical analysis, that this was not the case.
\item
When $\rho$ is finite, the critical covered volume looks concave as a function of $\alpha$.
However $\phi(\cdot,\infty)$ is not concave as soon as $\phi_c(\delta_1)<1-\exp(-2)$.
Based on \eqref{e:approx}, $\phi(\cdot,\infty)$ is therefore not concave.
As a consequence, at least for large enough $\rho$, $\phi(\cdot,\rho)$ is not concave.
\item
The numerical results suggests that the minimum of the critical covered fraction is reached when all the disks have the same radius.
(Equivalently, for all $\rho$ and all $\alpha$, $\phi(\alpha,\rho) \ge \phi(0,\rho)=\phi(1,\rho)=\phi_c(\delta_1)$.)
However there is neither a proof nor a disproof of such a result.
\item
The numerical results also suggest some monotonicity in $\rho$.
This has not been proven nor disproven.
\end{itemize}

\section{Proof of Theorem \ref{th:0}}
\label{s:preuve-th:0}

\subsection{Some notations}

In the whole of Section \ref{s:preuve-th:0}, we make the following assumptions:
\begin{itemize}
\item $\mu$ satisfies \eqref{e:th}.
\item $1<\widehat{\lambda_c}(\mu)$.
\end{itemize}

For all $\eta>0$,
we denote by $T_{\eta}\mu$ the measure defined by $T_{\eta}\mu(A)=\mu(A-\eta)$.
In other words, we can built $\Sigma(T_{\eta}\mu)$ from $\Sigma(\mu)$ by adding $\eta$ to each radius.

For all $\rho>1$,
we denote by $H^{\rho}\mu$ the measure defined by $H^{\rho}\mu(A)=\rho^d\mu(\rho A)$.
With this definition, $\rho^{-1}\Sigma(\mu)$ is a Boolean model driven by the measure $H^{\rho}\mu$.
For all $n \ge 0$, we let:
$$
m_n^{\rho} = \sum_{k=0}^n H^{\rho^k}\mu.
$$
With this definition and the notations of \eqref{e:definition-multi}, 
$$
\bigcup_{k=0}^n \rho^{-k} \Sigma_k
$$
is a Boolean model driven by $m_n^{\rho}$.
We also let:
$$
m_{\infty}^{\rho} = \sum_{k \ge 0} H^{\rho^k}\mu.
$$
So, $\Sigma^{\rho}(\mu)$ is a Boolean model driven by the locally finite measure $m_{\infty}^{\rho}$.

Let $p(a,\mu)$ denote the probability of existence of a path from $S(a/2)$ to $S(a)$ in $\Sigma(\mu)$:
$$
p(a,\mu) = P(\perco{S(a/2)}{S(a)}{\Sigma(\mu)}).
$$
We aim at proving that, for large enough $\rho$, $p(a,m_{\infty}^{\rho})\to 0$ as $a$ tends to infinity and $\Sigma^{\rho}(\mu)$ does not percolate. 
The first item of Theorem \ref{th:0} follows by applying this result to the measure $\lambda\mu$.
Recall that the second item of Theorem \ref{th:0} is contained in Theorem \ref{th:Gmulti}.

\subsection{Ideas}
\label{s:idees}

In this subsection we first sketch the proof of the existence of $\rho$ and $a$ such that $p(a,m_{\infty}^{\rho})$ is small.
This gives the main ingredients of the proof of the first item of Theorem \ref{th:0}.
A full proof is given in Subsection \ref{ss:proof}.
We then give the ideas of the proof of Theorems \ref{th:MPV1} and \ref{th:MPV2} by Menshikov, Popov and Vachkovskaia.
Their basic strategy is similar but the implementation of the proofs are different.

\subsubsection*{Sketch of the proof of the first item of Theorem \ref{th:0}}

Consider a small $\epsilon_1>0$.
Fix a small $\eta>0$ and a large $a$ such that (see Lemma \ref{l:eta}): 
\begin{equation}\label{e:maj0}
p(a,T_{\eta}\mu) \le \epsilon_1/2.
\end{equation}

For all $n \ge 1$, write:
$$
m_n^{\rho}=H^{\rho}m_{n-1}^{\rho}+\mu.
$$
If the event $\{\perco{S(a/2)}{S(a)}{\Sigma(m_n^{\rho})}\}$ occurs,
then either the event $\{\perco{S(a/2)}{S(a)}{T_{\eta}\mu}\}$ occurs (with a natural coupling between the Boolean models)
either in $\Sigma(H^{\rho}m_{n-1}^{\rho}) \cap B(a)$ one can find a component of diameter at least $\eta$.
We use this observation through its following crude consequence (see Lemma \ref{l:dichotomie}):
$$
p(a,m_n^{\rho}) \le p(a,T_{\eta}\mu) + Ca^d\eta^{-d}p(\eta/2,H^{\rho}m_{n-1}^{\rho}).
$$
By scaling and by \eqref{e:maj0}, this yields:
\begin{equation} \label{e:majn}
p(a,m_n^{\rho}) \le \epsilon_1/2 + Ca^d\eta^{-d}p(\rho\eta/2,m_{n-1}^{\rho}).
\end{equation}

But for any $\epsilon_2$, any small enough $\epsilon_1$ and any large enough $a$ we can find $\tau$ such that 
(see Lemmas \ref{l:carre} and \ref{l:poitiers}): 
\begin{equation}\label{e:turbo}
p(\tau a, m_{n-1}^{\rho}) \le \epsilon_2 \hbox{ as soon as } p(a, m_{n-1}^{\rho}) \le \epsilon_1.
\end{equation}
An important fact is that $\tau$ does not depends on $n$ nor on $\rho$,
provided $\rho \ge \rho_0$ where $\rho_0$ is an arbitrary constant strictly larger than $1$.
Here we use assumption \eqref{e:th} to bound error terms due to the existence of large balls.

We choose $\epsilon_2$ such that:
$$
Ca^d\eta^{-d}\epsilon_2 = \epsilon_1/2.
$$
We set $\rho=2\tau a/\eta$.
Then, \eqref{e:majn} and \eqref{e:turbo} can be rewritten as follows:
\begin{eqnarray}
p(a,m_n^{\rho})  & \le & \epsilon_1/2 + Ca^d\eta^{-d}p(\tau a,m_{n-1}^{\rho})  \label{e:majn2} \\
Ca^d\eta^{-d}p(\tau a,m_{n-1}^{\rho}) & \le  & \epsilon_1/2 \hbox{ as soon as } p(a, m_{n-1}^{\rho}) \le \epsilon_1 \label{e:turbo2}.
\end{eqnarray}
As moreover \eqref{e:maj0} implies $p(a,m_0^{\rho}) \le \epsilon_1$ we get, by induction and then sending $n$ to infinity
(see Lemma \ref{l:compacite}):
$$
p(a,m_{\infty}^{\rho}) \le \epsilon_1.
$$
The convergence of $p(a,m_{\infty}^{\rho})$ to $0$ is then extracted from the above result for a small enough $\epsilon_2$ and from
arguments behind \eqref{e:turbo} applied to $m_{\infty}^{\rho}$ and other $\epsilon$.

\subsubsection*{Sketch of the proofs of Theorems \ref{th:MPV1} and \ref{th:MPV2} by Menshikov, Popov and Vachkovskaia}
\label{s:MPV}
Let us quickly describe the ideas of the proofs of Menshikov, Popov and Vachkovskaia.
Those ideas are used in their papers \cite{Menshikov-al-multi} and \cite{Menshikov-al-multi-unbounded} through a discretization of space ; 
we describe them in a slightly more geometric way.
For simplicity we only consider two scales: $\rho^{-1}\Sigma_1$ and $\Sigma_0$.
For simplicity, we also assume that the radius is one in the unscaled model ($\mu=\lambda\delta_1$).
We assume that the scale factor $\rho$ is large enough.
Assume that $C$ is a connected component of $\rho^{-1}\Sigma_1 \cup \Sigma_0$ whose diameter is a least $\alpha$ 
(it can be much larger) for a small enough constant $\alpha>0$.
Then, $C$ is included in the union of the following kind of sets:
\begin{enumerate}
\item connected components of $\rho^{-1}\Sigma_1$ whose diameter is at least $\alpha$ ;
\item balls of $\Sigma_0$ enlarged by $\alpha$ (same centers but the radii are $1+\alpha$ instead of $1$).
\end{enumerate}
Then, they show that the union of all those sets is stochastically dominated by a Boolean model similar to $\Sigma_0$
but with radii enlarged by a factor $\alpha$ and with density of centers $1+\alpha'$ times the corresponding density for $\Sigma_0$
for a suitable $\alpha'>0$.
This part uses $\lambda<\widehat{\lambda_c}$.
In some sense, 
one can therefore control percolation in the union of two models by percolation in one model.
Iterating the argument with some care in the constants $\alpha$ and $\alpha'$, 
one sees that -- for large enough $\rho$ -- one can control percolation in the multiscale model by percolation in a subcritical model. 
This yields the result.



\subsection{Proof of Theorem \ref{th:0}}

\label{ss:proof}

As $1<\widehat{\lambda_c}(\mu)$, we know that $p(a,\mu)$ tends to $0$ as $a$ tends to infinity.
We need the following slightly stronger consequence.

\begin{lemma} \label{l:eta} There exists $\eta>0$ such that $p(a,T_{\eta}\mu)$ tends to $0$.
\end{lemma}
\proof Let $\epsilon>0$ and $x>0$.
We have:
\begin{eqnarray}
H^{1+\epsilon}T_{\epsilon^2}\mu([x,+\infty[)
 & = & (1+\epsilon)^d T_{\epsilon^2}\mu([x(1+\epsilon),+\infty[) \nonumber \\
 & = & (1+\epsilon)^d \mu([x(1+\epsilon)-\epsilon^2,+\infty[) \nonumber \\
 & \le & \kappa(\epsilon) (1+\epsilon)^d \mu([x,+\infty[) \label{e:domination}
\end{eqnarray}
where
$$
\kappa(\epsilon)=\frac{\mu(]0,+\infty[)}{\mu([\epsilon,+\infty[)}.
$$
The inequality is proven as follows. 
If $x \ge \epsilon$, then $[x(1+\epsilon)-\epsilon^2,+\infty[ \subset [x,+\infty[$ and the result follows from $\kappa(\epsilon) \ge 1$.
If, on the contrary, $x<\epsilon$,
then the left hand side is bounded above by $(1+\epsilon)^d \mu(]0,+\infty[)$ which is itself bounded above
by the right hand side.

Note that $\kappa(\epsilon)(1+\epsilon)^d$ tends to $1$ as $\epsilon$ tends to $0$.
Let us say that a measure $\nu$ is subcritical if $\widehat{\lambda_c}(\nu)>1$.
As $\mu$ is subcritical, we get that 
$\kappa(\epsilon)(1+\epsilon)^d \mu$ is subcritical for small enough $\epsilon$.
We fix such an $\epsilon$.
By \eqref{e:domination}  we can couple a Boolean model driven by $H^{1+\epsilon}T_{\epsilon^2}\mu$ 
and a Boolean model driven by $\kappa(\epsilon)(1+\epsilon)^d \mu$ in such a way that the first one is contained in the second one.
Therefore the first one is subcritical.
By scaling, a Boolean model driven by $T_{\epsilon^2}\mu$ is then subcritical. 
We take $\eta=\epsilon^2$. \finpreuve


\begin{lemma} \label{l:dichotomie} Let $\nu_1$ and $\nu_2$ be two finite measures on $]0,+\infty[$.
One has, for all $\eta>0$ and $a \ge 4\eta$:
$$
p(a,\nu_1+\nu_2) \le p(a,T_{\eta}\nu_1)+C_1 a^d\eta^{-d}p(\eta/2,\nu_2)
$$
where $C_1=C_1(d)>0$ depends only on the dimension $d$.
\end{lemma}

\proof
Let $(x_i)_{i\in I}$ be a family of points such that :
\begin{itemize}
\item The balls $B(x_i,\eta/4)$, $i\le I$, cover $B(a)$.
\item There are at most $C_1 a^d\eta^{-d}$ points in the family where $C_1=C_1(d)$ depends only on the dimension $d$.
\end{itemize}

We couple the different Boolean model as follows.
Let $\Sigma(\nu_1)$ be a Boolean model driven by $\nu_1$.
Let $\Sigma(\nu_2)$ be a Boolean model driven by $\nu_2$.
Assume that $\Sigma(\nu_1)$ and $\Sigma(\nu_2)$ are independent.
Then $\Sigma(\nu_1) \cup \Sigma(\nu_2)$ is a Boolean model driven by $\nu_1+\nu_2$.
We set $\Sigma(\nu_1 + \nu_2)=\Sigma(\nu_1) \cup \Sigma(\nu_2)$.
We also consider $\Sigma(T_{\eta}\nu_1)$, the Boolean model obtained by adding $\eta$ to the radius of each ball of $\Sigma(\nu_1)$.
Thus $\Sigma(T_{\eta}\nu_1)$ is driven by $T_{\eta}\nu_1$.

Let us prove the following property:
\begin{equation} \label{e:key-eta}
\{\perco{S(a/2)}{S(a)}{\Sigma(\nu_1+\nu_2)}\} 
\subset \{\perco{S(a/2)}{S(a)}{\Sigma(T_{\eta}\nu_1)}\} \cup \bigcup_{i\in I} \{\perco{S(x_i,\eta/4)}{S(x_i,\eta/2)}{\Sigma(\nu_2)}\}.
\end{equation}

Assume that $\Sigma(\nu_1+\nu_2)=\Sigma(\nu_1) \cup \Sigma(\nu_2)$ connects $S(a/2)$ with $S(a)$.
Recall $a \ge 4\eta$.
If the diameter of all connected components of $\Sigma(\nu_2) \cap B(a)$ are less or equal to $\eta$, 
then $\Sigma(T_{\eta}\nu_1)$ connects $S(a/2)$ with $S(a)$.
Otherwise, let $C$ be a connected component of $\Sigma(\nu_2) \cap B(a)$ with diameter at least $\eta$.
Let $x,y$ be two points of $C$ such that $\|x-y\|>\eta$.
The point $x$ belongs to a ball $B(x_i,\eta/4)$.
As $y$ does not belong to $B(x_i,\eta/2)$, 
the component $C$ connects $S(x_i,\eta/4)$ to $S(x_i,\eta/2)$.
Therefore, $\Sigma(\nu_2)$ connects $S(x_i,\eta/4)$ to $S(x_i,\eta/2)$.
We have proven \eqref{e:key-eta}. 
The lemma follows. \finpreuve

\medskip

The following lemma is essentially the first item of Proposition 3.1 in \cite{G-perco-boolean-model}.
For the sake of completeness we nevertheless provide a proof.

\begin{lemma} \label{l:carre}
Let $\nu$ be a finite measure on $]0,+\infty[$.
There exists a constant $C_2=C_2(d)>0$ such that, for all $a>0$:
$$
p(10a,\nu) \le C_2 p(a, \nu)^2+ C_2 \int_{[a,+\infty[} r^d \nu(dr).
$$
\end{lemma}
\proof
Let $K$ be a finite subset of $S(5)$ such that $K+B(1/2)$ covers $S(5)$.
Let $L$ be a finite subset of $S(10)$ such that $L+B(1/2)$ covers $S(10)$. 
Let $A$ be the following event: there exists a random ball $B(c,r)$ of $\Sigma(\nu)$ such that $r\ge a$ and $B(c,r)\cap B(10a)$ is non empty.
We have:
$$
\{\perco{S(5a)}{S(10a)}{\Sigma(\nu)}\} \setminus A \subset \{\percorayon{S(5a)}{S(10a)}{\Sigma(\nu)}{a}\}
$$
where, in the last event, we ask for a path using only balls of $\Sigma(\nu)$ of radius at most $a$.
Let us prove the following:
\begin{multline}\label{e:key-carre}
\{\perco{S(5a)}{S(10a)}{\Sigma(\nu)}\} \setminus A \\
\subset
\bigcup_{k \in K, l \in L} 
\{\percorayon{S(ak,a/2)}{S(ak,a)}{\Sigma(\nu)}{a}\} \cap \{\percorayon{S(al,a/2)}{S(al,a)}{\Sigma(\nu)}{a}\}.
\end{multline}
Assume that the event on the left hand side occurs.
Then, by the previous remark, there exists a path from a point $x\in S(5a)$ to a point $y\in S(10a)$ 
that is contained in balls of $\Sigma(\nu)$ of radius at most $a$.
As $Ka+B(a/2)$ covers $S(5a)$, there exists $k\in K$ such that $x$ belongs to $B(ka,a/2)$.
Using the previous path, one gets that the event 
$$
\{\percorayon{S(ak,a/2)}{S(ak,a)}{\Sigma(\nu)}{a}\}
$$ 
occurs.
By a similar arguments involving $y$ we get \eqref{e:key-carre}.

Observe that, for all $k \in K$ and $l\in L$, the events 
$$
\{\percorayon{S(ak,a/2)}{S(ak,a)}{\Sigma(\nu)}{a}\} \hbox{ and } 
\{\percorayon{S(al,a/2)}{S(al,a)}{\Sigma(\nu)}{a}\}
$$
are independent.
Indeed, the first one depends only on balls with centers in $B(ak,2a)$, the second one depends only on balls with centers in $B(al,2a)$, and
$\|ak-al\| \ge 5a$.
Using this independence, stationarity and \eqref{e:key-carre}, we then get:
$$
P(\{\perco{S(5a)}{S(10a)}{\Sigma(\nu)}\}) \le CP(\percorayon{S(a/2)}{S(a)}{\Sigma(\nu)}{a}\})^2+P(A)
$$
where $C$ is the product of the cardinality of $K$ by the cardinality of $L$.
The probability $P(A)$ is bounded above by standard computations. \finpreuve

\medskip

From the previous lemma, we deduce the following result.

\begin{lemma} \label{l:poitiers} 
Let $\epsilon>0$. 
There exists $C_3=C_3(d)>0$, $a_0=a_0(d,\mu)$ and $k_0=k_0(d,\mu,\epsilon)$ such that, for all $N$, all $\rho \ge 2$ and all $a \ge a_0$:
if $p(a,m^{\rho}_N) \le C_3$ then for all $k \ge k_0$, $p(a10^k,m^{\rho}_N) \le \epsilon$.
\end{lemma}
\proof 
For all $\rho \ge 2$ and all $a \ge 1$ we have: 
\begin{eqnarray*}
\int_{[a,+\infty[} r^d m^{\rho}_{\infty}(dr) 
 & =  & \sum_{k\ge 0} \rho^{kd} \int_{]0,+\infty[} 1_{[a,+\infty[}(r \rho^{-k}) (r \rho^{-k})^d \mu(dr) \\
 & =  & \int_{]0,+\infty[} \sum_{k\ge 0} 1_{[a,+\infty[}(r \rho^{-k}) r^d \mu(dr) \\
 & =  & \int_{[a,+\infty[} \big(\left\lfloor \ln(r/a)\ln(\rho)^{-1} \right\rfloor +1\big) r^d \mu(dr) \\
 & \le & \int_{[a,+\infty[} (\ln(r)\ln(2)^{-1}+1) r^d \mu(dr).
\end{eqnarray*}

Let $C_2$ be the constant given by Lemma \ref{l:carre}.
By \eqref{e:th} we can chose $a_0=a_0(d,\mu) \ge 1$ such that 
\begin{equation}\label{e:defa0}
C_2^2 \int_{[a_0,+\infty[} (\ln(r)\ln(2)^{-1}+1) r^d \mu(dr) \le \frac{1}{4}.
\end{equation}
Let $C_3=(2C_2)^{-1}$. 
Let $N$, $\rho$ and $a$ be as in the statement of the lemma.
From Lemma \ref{l:carre} we get: 
\begin{eqnarray}
C_2 p(10a,m^{\rho}_N) 
 & \le & (C_2 p(a,m^{\rho}_N))^2+C_2^2 \int_{[a,+\infty[} r^d m^{\rho}_N(dr) \\
 & \le & (C_2 p(a,m^{\rho}_N))^2+C_2^2 \int_{[a,+\infty[} (\ln(r)\ln(2)^{-1}+1) r^d \mu(dr) \label{e:angouleme}
\end{eqnarray}
Let $(u_k)$ be a sequence defined by $u_0=1/2$ and, for all $k\ge 0$:
\begin{equation}\label{e:defu}
u_{k+1} = u_k^2+C_2^2 \int_{[a_0 10^k,+\infty[} (\ln(r)\ln(2)^{-1}+1) r^d \mu(dr).
\end{equation}
Note that the sequence $(u_k)$ only depends on $d$ and $\mu$.

Assume that $p(a,m^{\rho}_N)\le C_3$. 
We then have $C_2p(a,m^{\rho}_N) \le u_0$.
Using $a\ge a_0$ and \eqref{e:angouleme}, we then get $C_2p(a10^k,m^{\rho}_N) \le u_k$ for all $k$.
Therefore, it sufficies to show that the sequence $(u_k)$ tends to $0$.

Using \eqref{e:defu}, \eqref{e:defa0} and $u_0=1/2$ we get $0 \le u_k \le 1/2$ for all $k$. 
Therefore, $0 \le \limsup u_k \le 1/2$.
By \eqref{e:defu} and by the convergence of the integrale we also get $\limsup u_k \le (\limsup u_k)^2$.
As a consequence, $\limsup u_k=0$ and the lemma is proven. \finpreuve

\begin{lemma} \label{l:compacite} For all $a>0$ and $\rho>1$ the following convergence holds:
$$
p(a,m^{\rho}_{\infty})=\lim_{N\to\infty} p(a,m^{\rho}_N) .
$$
\end{lemma}
\proof 
The sequence of events 
$$
A_N=\{\perco{S(a/2)}{S(a)}{\Sigma(m^{\rho}_N)}\}
$$
is increasing (we use the natural coupling between our Boolean models). 
Therefore, it suffices to show that the union of the previous events is 
$$
A=\{\perco{S(a/2)}{S(a)}{\Sigma(m^{\rho}_{\infty})}\}.
$$
If $A$ occurs, then there is is path from $S(a/2)$ to $S(a)$ that is contained in $\Sigma(m^{\rho}_{\infty})$.
By a compactness argument, this path is included in a finite union of ball of $\Sigma(m^{\rho}_{\infty})$.
Therefore, there exists $N$ such that the path is included in $\Sigma(m^{\rho}_N)$ and $A_N$ occurs.
This proves $A \subset \cup A_N$. 
The other inclusion is straightforward. \finpreuve

\proofof{the second item of Theorem \ref{th:0}}
By Lemma \ref{l:eta}, we can fix $\eta_1>0$ such that $p(a,T_{10 \eta_1}\mu)$ tends to $0$ as $a$ tends to $\infty$.
Let $C_1$ be given by Lemma \ref{l:dichotomie}.
Let $a_0$ and $C_3$ be as given by Lemma \ref{l:poitiers}.
Fix $a_1 \ge \max(40\eta_1,a_0,1)$  such that $p(a,T_{10\eta_1}\mu) \le C_3/2$ for all $a \ge a_1$.
Let $k_0$ be given by Lemma \ref{l:poitiers} with the choice:
$$
\epsilon=C_1^{-1}(10a_1)^{-d}\eta_1^dC_3/2.
$$
Therefore, for all $\rho \ge 2$, all $N$, all $a\in [a_1, 10a_1]$ and all $\eta\in [\eta_1, 10\eta_1]$: 
$$
C_1 a^d \eta^{-d} p(a 10^k,m_N^{\rho}) \le \frac{C_3}{2} \hbox{ for all } k \ge k_0 \hbox{ as soon as } p(a,m_N^{\rho}) \le C_3.
$$
Fix $k \ge k_0$, $a \in [a_1, 10a_1]$ and $\eta \in [\eta_1, 10\eta_1]$.
Set:
$$
\rho=2a 10^k \eta^{-1}.
$$
Note $\rho \ge 8 \ge 2$ as $a \ge a_1 \ge 40\eta_1 \ge 4\eta$.
By Lemma \ref{l:dichotomie} we have, for all $N$:
\begin{eqnarray*}
p(a,m_{N+1}^{\rho}) 
 & \le & p(a,T_{\eta}\mu) + C_1 a^d \eta^{-d} p(\eta/2, m^{\rho}_{N+1}-\mu) \\
 & =  & p(a,T_{\eta}\mu) + C_1 a^d \eta^{-d} p(\eta/2, H^{\rho} m_N^{\rho}).
\end{eqnarray*}
By definition of $a_1$, by $a_1 \le a$, by $\eta \le 10\eta_1$, by scaling and by definition of $\rho$ we get, for all $N$:
\begin{eqnarray*}
p(a,m_{N+1}^{\rho}) 
 & \le & \frac{C_3}{2}+C_1 a^d \eta^{-d} p(\rho\eta/2, m_N^{\rho}) \\
 & = & \frac{C_3}{2}+C_1 a^d \eta^{-d} p(a 10^k, m_N^{\rho}).
\end{eqnarray*}
Combining this inequality with the property defining $k_0$, we get that
 $p(a,m_N^{\rho}) \le C_3$ implies $p(a,m_{N+1}^{\rho}) \le C_3$.
As $p(a,m_0^{\rho}) = p(a,\mu) \le p(a,T_{10\eta_1}\mu) \le C_3/2$ we get $p(a,m_N^{\rho}) \le C_3$ for all integer $N$.

Let $\epsilon>0$. 
Using again Lemma \ref{l:poitiers} we get the existence of an integer $k'_0$ such that 
$p(a 10^{k'}, m_N^{\rho}) \le \epsilon$ for all $k' \ge k'_0$ as soon as $p(a,m_N^{\rho}) \le C_3$. 
But we have proven the latter property. 
Therefore $p(a 10^{k'}, m_N^{\rho}) \le \epsilon$ for all $N$ and all $k' \ge k'_0$.
By Lemma \ref{l:compacite}, we get $p(a 10^{k'}, m_{\infty}^{\rho}) \le \epsilon$ for all $k' \ge k'_0$.
Using the freedom on the choice of $k \ge k_0$ and $\eta \in [\eta_1, 10\eta_1]$, we get that the previous result holds for all 
$\rho \ge 2a 10^{k_0-1} \eta_1^{-1}$ and then for all $\rho \ge 2a_1 10^{k_0} \eta_1^{-1}$.
Moreover, using the freedom on the choice of $a \in [a_1, 10a_1]$ and $k' \ge k'_0$, we get:
$$
p(r,m_{\infty}^{\rho}) \le \epsilon \hbox{ for all } r \ge a_1 10^{k'_0} \hbox{ and all } \rho\ge 2a_1 10^{k_0} \eta_1^{-1}.
$$
Therefore, $p(r,m_{\infty}^{\rho})$ tends to $0$ as $r$ tends to infinity.
As a consequence, $\Sigma^{\rho}(\mu)$ does not percolate for any $\rho \ge 2a_1 10^{k_0} \eta_1^{-1}$.
\finpreuve

\section{Proof of Theorem \ref{th:1}}
\label{s:preuve-th:1}

\begin{lemma} \label{l:s}
Let $s>0$ and $\rho>1$.
The following assumptions are equivalent:
\begin{enumerate}
\item $\int_{]0,+\infty[} r^{d+s} \mu(dr) < \infty$.
\item $\int_{[1,+\infty[} r^{d+s} m_{\infty}^{\rho}(dr) < \infty$.
\end{enumerate}
\end{lemma}
\proof
We have:
\begin{eqnarray*}
\int_{[1,+\infty[} r^{d+s} m^{\rho}_{\infty}(dr) 
 & =  & \sum_{k\ge 0} \rho^{kd} \int_{]0,+\infty[} 1_{[1,+\infty[}(r \rho^{-k}) (r \rho^{-k})^{d+s} \mu(dr) \\
 & =  & \int_{[1,+\infty[} \sum_{k\ge 0} 1_{[1,+\infty[}(r \rho^{-k}) \rho^{-ks} r^{d+s} \mu(dr).
\end{eqnarray*}
Therefore:
$$
\int_{[1,+\infty[} r^{d+s} \mu(dr) \le  \int_{[1,+\infty[} r^{d+s} m^{\rho}_{\infty}(dr) \le \frac{1}{1-\rho^{-s}} \int_{[1,+\infty[} r^{d+s} \mu(dr).
$$
This yields the result.
\finpreuve

\proofof{the first item of Theorem \ref{th:1}}
By the discussion at the beginning of Section 1.5 in \cite{G-perco-generale},
$\Sigma^{\rho}(\lambda\mu)$ is driven by a a Poisson point process whose intensity is the product
of the Lebesgue measure by the locally finite measure $\lambda m_{\infty}^{\rho}$.
Let us check the three items of Theorem 2.9 in \cite{G-perco-generale} with $\rho=10$
($\rho$ is not use in the same way in \cite{G-perco-generale}).
We refer to Section 2.1 of \cite{G-perco-generale} for definitions.
\begin{enumerate}
\item 
The first item is fulfilled thanks to \eqref{e:pitilde0} 
\item 
For all $\beta>0$ and all $x \in \R^d$, 
the event $G(x,0,\beta)$ only depends on balls $B(c,r) \in \Sigma^{\rho}(\lambda\mu)$ such that $c$ belongs to $B(x,3\beta)$.  
By the independance property of Poisson point processes, we then get that $G(0,0,\beta)$ and $G(x,0,\beta)$ are independent whenever $\|x\| \ge 10\beta$.
Therefore $I(10,0,\beta)=0$ and the second item of Theorem 2.9 is fulfilled.
\item 
The third item (note that $\mu$ in \cite{G-perco-generale} is $m_{\infty}^{\rho}$ in this paper) is fulfilled thanks to Lemma \ref{l:s}
\end{enumerate}
Theorem 2.9 in \cite{G-perco-generale} yields the result. \finpreuve

\proofof{the second item of Theorem \ref{th:1}}
If $\int r^d\mu(dr)$ is infinite then, $\Sigma(\lambda\mu)$ percolates for all $\lambda>0$ (see the dicussion of Section \ref{s:boolean}). 
Therefore $\Sigma^{\rho}(\lambda\mu)$ percolates for all $\rho>1$ and $\lambda>0$.
Therefore $D^{\rho}(\lambda\mu)=\infty$ with positive probability for all $\rho>1$ and $\lambda>0$.  

Now, assume that $\int r^d \mu(dr)$ is finite. 
Then, by the discussion at the beginning of Section 1.5 in \cite{G-perco-generale},
$\Sigma^{\rho}(\lambda\mu)$ is driven by a a Poisson point process whose intensity is the product
of the Lebesgue measure by the locally finite measure $\lambda m_{\infty}^{\rho}$.
We can therefore apply Theorem 1.2 in \cite{G-perco-generale}.
By Lemma \ref{l:s}, assumption (A3) of Theorem 1.2 in \cite{G-perco-generale} is not fulfilled
(note that $\mu$ in \cite{G-perco-generale} is $m_{\infty}^{\rho}$ in this paper).
Theorem 1.2 in \cite{G-perco-generale} then yields the result. \finpreuve

\section{Proof of Proposition \ref{p}}
\label{s:preuve-p}

We first need a lemma, which is a consequence of Lemmas \ref{l:dichotomie} and \ref{l:carre}.

\begin{lemma} \label{l:resume} 
Let $\nu_1$ and $\nu_2$ be two finite measures on $]0,+\infty[$.
Let $\eta>0$ and $a_0 \ge 4\eta$.
Let $\rho>1$.
There exists $C_4=C_4(d)>0$ such that $\widehat{\lambda}_c(\nu_1+H_{\rho}(\nu_2)) \ge 1$ as soon as the following conditions hold:
\begin{enumerate}
\item $p(a,T_{\eta}\nu_1) \le C_4$ for all $a \in [a_0, 10a_0]$.
\item $a_0^d \eta^{-d} p(\rho\eta/2, \nu_2) \le C_4$.
\item $\int_{[a_0,+\infty[} r^d \nu_1(dr) \le C_4$ and $\int_{[a_0,+\infty[} r^d \nu_2(dr) \le C_4$.
\end{enumerate}
\end{lemma}

\proof Let $C_4=C_4(d)>0$ be such that $C_4C_2(1+10^dC_1) \le 1/2$ and $2C_2^2C_4 \le 1/4$,
where $C_1$ appears in Lemma \ref{l:dichotomie} and $C_2$ appears in Lemma \ref{l:carre}.
Set $\nu=\nu_1+H_{\rho}(\nu_2)$. 

For all $a \in [a_0,10a_0]$ we have, 
by Lemma \ref{l:dichotomie} applied to $\nu_1$ and $H_{\rho}(\nu_2)$, by scaling and by the assumptions of the lemma:
\begin{eqnarray}
p(a,\nu) 
 & \le & p(a,T_{\eta}\nu_1) + C_1 a^d\eta^{-d}p(\eta/2,H_{\rho}\nu_2) \nonumber \\
 & = & p(a,T_{\eta}\nu_1) + C_1 a^d\eta^{-d}p(\rho\eta/2,\nu_2) \nonumber \\ 
 & \le & C_4(1+10^dC_1) \nonumber \\
 & \le & 1/(2C_2). \label{e:fevrier0}
\end{eqnarray}

But for all $a \ge a_0$ we have, by Lemma \ref{l:carre} and by the assumptions of the lemma:
\begin{eqnarray}
C_2p(10a,\nu) 
 & \le & (C_2p(a,\nu))^2 + C_2^2\int_{[a,+\infty[}r^d \nu(dr) \nonumber \\
 & = & (C_2p(a,\nu))^2 + C_2^2\int_{[a,+\infty[}r^d \nu_1(dr) + C_2^2\int_{[a\rho,+\infty[}r^d \nu_2(dr) \label{e:fevrier1} \\ 
 & \le & (C_2p(a,\nu))^2 + 2C_2^2C_4 \nonumber \\
 & \le & (C_2p(a,\nu))^2 + 1/4.\label{e:fevrier2}
\end{eqnarray}

By \eqref{e:fevrier0} and \eqref{e:fevrier2} we get $C_2p(a,\nu) \le 1/2$ for all $a \ge a_0$ and therefore $0 \le \limsup C_2p(a,\nu) \le 1/2$.
By \eqref{e:fevrier1} and the third assumption of the lemma, we get $\limsup C_2p(a,\nu) \le (\limsup C_2p(a,\nu))^2$.
Therefore, we must have $\limsup C_2p(a,\nu)=0$ 
and the lemma is proven. \finpreuve

\proofof{Proposition \ref{p}}
The inequality is straightforward.
To prove the inequality, we note that, by scaling, $\widehat{\lambda}_c(H_{\rho}\nu_2)=\widehat{\lambda}_c(\nu_2)$.
Let us prove the convergence.
We can assume $\widehat{\lambda}_c(\nu_1)>0$ and $\widehat{\lambda}_c(\nu_2)>0$,
otherwise the convergence is obvious.
Therefore, by Lemma \ref{l:hat}, the integrals $\int r^d \nu_1(dr)$ and $\int r^d \nu_2(dr)$ are finite.

Let $C_4$ be the constant given by Lemma \ref{l:resume}.
Let $0<\epsilon<1$.
Note: 
$$
\widehat{\lambda_c}\big((1-\epsilon)\widehat{\lambda}_c(\nu_1)\nu_1\big) = (1-\epsilon)^{-1}\widehat{\lambda}_c(\nu_1)^{-1}\widehat{\lambda}_c(\nu_1)>1,
$$
Therefore, by Lemma \ref{l:eta} (in which \eqref{e:th} is not used), we can fix $\eta>0$ such that
$$
p(a,T_{\eta}(1-\epsilon)\widehat{\lambda}_c(\nu_1)\nu_1) \to 0.
$$
We can then fix $a_0 \ge 4\eta$ such that:
\begin{equation}\label{e:annexe1}
p(a,T_{\eta}(1-\epsilon)\widehat{\lambda}_c(\nu_1)\nu_1) \le C_4 \hbox{ for all } a \ge a_0
\end{equation}
and such that 
\begin{equation}\label{e:annexe3.1}
\int_{[a_0,+\infty[} r^d \widehat{\lambda}_c(\nu_1)\nu_1(dr) \le C_4
\end{equation}
and
\begin{equation}\label{e:annexe3.2}
\int_{[a_0,+\infty[} r^d \widehat{\lambda}_c(\nu_2)\nu_2(dr) \le C_4.
\end{equation}
Now we fix $\rho_0>1$ such that :
\begin{equation}\label{e:annexe2}
a_0^d\eta^{-d} p(\rho\eta/2, (1-\epsilon)\widehat{\lambda}_c(\nu_2)\nu_2) \le C_4 \hbox{ for all } \rho\ge \rho_0.
\end{equation}
Now, let $0<\alpha<1$ and let 
$$
\lambda = \min\left(\frac{\widehat{\lambda}_c(\nu_1)(1-\epsilon)}{\alpha}, \frac{\widehat{\lambda}_c(\nu_2)(1-\epsilon)}{1-\alpha}\right).
$$
By \eqref{e:annexe1}, \eqref{e:annexe2}, \eqref{e:annexe3.1} and \eqref{e:annexe3.2} 
we get that Assumptions 1 
, 2 and 3 of Lemma \ref{l:resume} are fulfilled
for the measures $\alpha\lambda\nu_1$ and $(1-\alpha) \lambda \nu_2$ and for $\rho \ge \rho_0$.
Therefore, we get 
$$
\widehat{\lambda}_c(\alpha \lambda \nu_1 + (1-\alpha) \lambda H_{\rho} \nu_2) \ge 1
$$
and thus:
$$
\widehat{\lambda}_c(\alpha \nu_1 + (1-\alpha)  H_{\rho} \nu_2) \ge \lambda = (1-\epsilon)
\min\left(\frac{\widehat{\lambda}_c(\nu_1)}{\alpha}, \frac{\widehat{\lambda}_c(\nu_2)}{1-\alpha}\right).
$$
Therefore, as soon as $\rho \ge \rho_0$, we have:
\begin{eqnarray*}
0 
& \le &  
\min\left(\frac{\widehat{\lambda}_c(\nu_1)}{\alpha}, \frac{\widehat{\lambda}_c(\nu_2)}{1-\alpha}\right) 
- \widehat{\lambda}_c(\alpha \nu_1 + (1-\alpha)  H_{\rho} \nu_2) \\
& \le & 
\epsilon \min\left(\frac{\widehat{\lambda}_c(\nu_1)}{\alpha}, \frac{\widehat{\lambda}_c(\nu_2)}{1-\alpha}\right)  \\
& \le &
\epsilon \max(2\widehat{\lambda}_c(\nu_1),2\widehat{\lambda}_c(\nu_2)).
\end{eqnarray*}
This yields the proposition. \finpreuve

\appendix

\section{Critical parameters}

\begin{lemma} \label{l:lambdacinegalite}
$$
\widetilde{\lambda_c}(\mu) \le \widehat{\lambda_c}(\mu) \le \lambda_c(\mu). 
$$
\end{lemma}
\proof
The second inequality is a consequence of the following inclusion:
$$
\{\perco{\{0\}}{S(r)}{\Sigma}\} \subset \{\perco{S(r/2)}{S(r)}{\Sigma}\}.
$$
The first inequality can be proven as follows.
Let $r \ge 1$.
By the FKG inequality, we get:
\begin{eqnarray*}
P(\perco{\{0\}}{S(r)}{\Sigma}) 
 & \ge & P(B(0,1) \subset \Sigma \hbox{ and } \perco{S(1)}{S(r)}{\Sigma} ) \\
 & \ge & C P(\perco{S(1)}{S(r)}{\Sigma} )
\end{eqnarray*}
where $C = P(B(0,1) \subset \Sigma)>0$ does not depend on $r$.
For all large enough $r$, we can cover $S(2r)$  by at most $C'r^d$ balls $B(x_i,1)$ where $C'$ only depends on the dimension $d$.
If there is a path in $\Sigma$ from $S(2r)$ to $S(4r)$, 
then there exists $i$ and a path in $\Sigma$ from $S(x_i,1)$ to $S(x_i,r)$.
(Consider the ball $B(x_i,1)$ that contains the initial point of the path.)
By stationarity and by the previous inequality we thus get:
\begin{eqnarray*}
P(\perco{S(2r)}{S(4r)}{\Sigma}) 
 & \le & C'r^d P(\perco{S(1)}{S(r)}{\Sigma}) \\
 & \le & C'C^{-1}r^dP(\perco{\{0\}}{S(r)}{\Sigma}).
\end{eqnarray*}
The first inequality stated in the lemma follows. \finpreuve

\begin{lemma} \label{l:hat} The treshold parameter $\widehat{\lambda_c}(\mu)$ is positive if and only if $\int r^d \mu(dr)$ is finite.
\end{lemma}
\proof If $\widehat{\lambda_c}(\mu)$ is positive, then there exists $\lambda>0$ such that $\Sigma(\lambda\mu)$ does not percolate.
By Theorem 2.1 of \cite{G-perco-boolean-model} this implies that $\int r^d \mu(dr)$ is finite.

Let us assume now that $\int r^d \mu(dr)$ is finite.
We need to prove the existence of $\lambda>0$ such that $p(a,\lambda\mu)$ tends to $0$.
This is proven, as an intermediate result, in the proof of Theorem 1.1 in \cite{G-perco-generale}.
As the result is an easy consequence of Lemma \ref{l:carre}, we find it more convenient to provide a proof here.
Let $C_2$ be the constant given by Lemma \ref{l:carre}.
For all $a>0$ and $\lambda>0$ we have:
\begin{equation}\label{e:flo1}
C_2p(10a,\lambda\mu) \le (C_2p(a,\lambda\mu))^2+\lambda C_2^2\int_{[a,+\infty[} r^d\mu(dr).
\end{equation}
For all $0<a \le 1$ we have, by standard computations:
$$
C_2p(a,\lambda\mu) \le C_2 P(\hbox{a ball of }\Sigma(\lambda\mu)\hbox{ touches }B(a)) \le C_2v_d \lambda \int_{]0,+\infty[} (1+r)^d \mu(dr)  
$$
where $v_d$ is the volume of the unit Euclidean ball.
As $\int r^d \mu(dr)$ is finite, we can therefore fix $\lambda>0$ such that:
\begin{equation}\label{e:flo2}
\lambda C_2^2\int_{[a,+\infty[} r^d\mu(dr) \le 1/4 \hbox{ for all } a>0 \hbox{ and } C_2p(a,\lambda\mu) \le 1/2 \hbox{ for all } 0<a \le 1.
\end{equation}
By \eqref{e:flo1}, \eqref{e:flo2} and by induction we get $C_2p(a,\lambda\mu) \le 1/2$ for all $a>0$. 
Therefore, we have $0 \le \limsup C_2p(a,\lambda\mu) \le 1/2$.
But \eqref{e:flo1} also yields the inequality $\limsup C_2p(a,\lambda\mu) \le (\limsup C_2p(a,\lambda\mu))^2$.
As a consequence we must have $\limsup C_2p(a,\lambda\mu)=0$ and then $p(a,\lambda\mu) \to 0$. \finpreuve


\def\cprime{$'$} \def\cprime{$'$}

\end{document}